\documentclass[10pt,a4paper]{article}

\usepackage{pifont,latexsym,amssymb,bbm}
\usepackage{yfonts,url}

\usepackage{hyperref}

\usepackage{graphicx}

\newtheorem{theorem}{Theorem}

\def\bd{\begin{description}}
\def\ed{\end{description}}

\def\beq{\begin{equation}}
\def\eeq{\end{equation}}
\def\bea{\begin{eqnarray}}
\def\eea{\end{eqnarray}}
\def\beas{\begin{eqnarray*}}
\def\eeas{\end{eqnarray*}}

\def\G1{\hbox{$\displaystyle{\mbox{\ding{172}}}$}}

\begin{document}

\title{Lower and upper estimates of the quantity\\ of algebraic numbers}

\author{  Yaroslav D. Sergeyev\\ University of Calabria,  Rende, Italy  \\   Lobachevsky University,
    Nizhni Novgorod, Russia   \\  Institute of High Performance Computing and
  Networking\\ of the National Research Council of Italy, Rende, Italy \\
  \\ https://www.yaroslavsergeyev.com
}
 \date{}

\maketitle


\begin{abstract}It is well known   that the set   of algebraic numbers (let us call it $A$)  is countable. In this paper, instead of the usage of the classical terminology of cardinals proposed by Cantor, a recently introduced   methodology using \G1-based infinite numbers is applied to measure the set $A$ (where the number \G1 is called \emph{grossone}). Our interest to this methodology is explained by the fact that in certain   cases where cardinals allow one to say only whether a set is countable or it has the cardinality of the continuum, the \G1-based  methodology   can provide a more accurate measurement  of infinite sets. In this article, lower and upper estimates of the number of elements of $A$ are obtained. Both estimates are expressed in \G1-based numbers.
\end{abstract}



\section{Introduction}

The notion of algebraic numbers is very well known: a complex number~$z$ is called \emph{algebraic} if there exist integers $a_i \in \mathbb{Z}, i = 0, \ldots n$, not all equal to zero, such that
\[
a_0 z^n + a_1 z^{n-1} + \ldots +  a_{n-1}  z + a_n =0.
\]
Let us call the set of all algebraic numbers $A$. It is easy to show that this set is infinite. In fact, it is sufficient to consider the roots of the polynomial $z-n=0$, $n \in \mathbb{Z}$. Since its unique solution is $z=n$, it follows that $\mathbb{Z} \subset A$. Insofar as the set of integers $\mathbb{Z}$ is infinite, the set $A$ is also infinite. Moreover, it is also known   that the set  $A$   is countable, i.e., its cardinality is $\aleph_0$. One of possible proofs can be sketched as follows.

Let       $B_{k}$ be the set of all  tuples of integers
\beq
(a_0, a_1, \ldots, a_k), \,\, a_i \in \mathbb{Z}, \, i=0,\ldots, k,\,\,a_0 \neq 0,  \label{tuples_1}
\eeq
  and the elements $a_0, a_1, ...,\ldots, a_k$ do not need to be distinct. Then, clearly, the set $B_{k}$ is countable.  For each tuple $(a_0, a_1, \ldots, a_k)  \in B_k$ let us consider  the corresponding polynomial
\beq
a_0 z^k + a_1 z^{k-1} + \ldots +  a_{k-1}  z + a_k =0.
\label{polyn}
\eeq
From the fundamental theorem of algebra, we know that there are exactly~$k$ complex roots (counted according to their multiplicities) for each polynomial (\ref{polyn}).

Since $k$ is a natural number and the set of natural numbers is countable,  we have a countable number of sets $B_{k}$, each containing a countable number of tuples (\ref{tuples_1}), each of which corresponds to $k$ roots of a $k$-degree polynomial.   The set $A$ is, therefore, a countable union of a countable  union of roots of polynomials corresponding to tuples (\ref{tuples_1}) of sets $B_{k}$, i.e., $A$ is also countable.

In this paper, instead of the usage of the classical terminology of cardinals   (see
\cite{Cantor} and modern developments in \cite{Bagaria_1,Heller,Jech,Mancosu_2016,Ternullo_Fano}), to measure the set $A$ we apply a recently introduced counting methodology (see a comprehensive survey in
\cite{EMS} and a  popular presentation  in \cite{Sergeyev}). Our interest to this methodology is explained by the fact that there exist cases where cardinals allow one to say only whether certain sets are countable or they have the cardinality of the continuum. In its turn, the new methodology using an infinite unit of measure called \emph{grossone} and expressed by the symbol~\G1  can provide   a more accurate measurement  of these infinite sets.

Hereinafter we suppose that the reader is familiar with the \G1-based methodology. Only a list of some applications of this methodology and  a few notions and results from \cite{EMS} strictly necessary for the further  presentation  will be recalled in Section~2. The main results providing lower and upper estimates for the quantity   of algebraic numbers are presented in Section~3. Finally, Section~4 concludes the paper.

\section{A theoretical background}

 The \G1-based methodology    has attracted a lot of attention from scientists working in different areas of mathematics and computer science. We provide here just a few examples of areas where this methodology is useful. First of all, we mention numerous applications in  local, global, and
multi-criteria optimization and classification   (see, e.g.,
\cite{Astorino,lexBB,Cococcioni2021big,DeLeone2020b} and references given therein). Then, we can indicate game theory (see, e.g., \cite{DAlotto_games,Fiaschi2020}),  probability theory (see, e.g., \cite{Calude_Dumitrescu,Pepelyshev_Zhigljavsky,Rizza_2,Rizza_3}),  fractals (see, e.g.,
\cite{Antoniotti2020,Caldarola2020_Soft_Comp,Koch}),  infinite
series   (see
\cite{EMS,paradoxes,Zhigljavsky}),
  Turing machines, cellular automata, and ordering  (see, e.g.,
\cite{DAlotto,Rizza_3,medals,Sergeyev_Garro_2}),   numerical
differentiation and numerical solution of  ordinary differential
equations (see, e.g.,  \cite{ODE_3,Falcone:et:al.2021,Falcone2022,ODE_5} and references given therein), etc.

In particular,   successful applications of this methodology in teaching mathematics should be mentioned   (see  \cite{Antoniotti_teaching,Iannone,Ingarozza_teaching,Mazzia_teaching,Nasr_1,Sergeyev_teaching}). The dedicated web page \cite{web_teaching}, developed at the University of East Anglia, UK, contains, among other things, a comprehensive teaching manual. It should be also emphasized that numerous papers studying consistency of the new methodology  and its connections   to the historical panorama of ideas dealing with infinities and
infinitesimals have been
published (see
\cite{Gangle,Lolli_2,MM_bijection,Sorbi,first,Fallacies,Tohme_univalent}). In particular, it is stressed in \cite{Fallacies} that this methodology is not related to non-standard analysis. Many other applications from different fields of computer science and pure and applied mathematics can be found at the web page \cite{www}.

 \begin{table}[t!]
 \caption{Cardinalities and the \G1-based numbers of elements of some infinite sets (see \cite{EMS}, p. 287).}
\begin{center} \small \label{table1}
\begin{tabular}{@{\extracolsep{\fill}}|c|c|c| }\hline
  Description     &  Cantor's  & Number  of \\
 of sets  &   cardinalities  &   elements  \\
 \hline
  &   &     \vspace{-2mm}   \\
the set of natural numbers  $\mathbb{N}$   &  countable, $\aleph_0$ & \G1     \\
   &   &     \vspace{-2mm}   \\
  $\mathbb{N} \setminus \{ 3, 5, 10, 23, 114 \} $   &  countable, $\aleph_0$ & \G1-5     \\
  &   &     \vspace{-2mm}   \\
 the set of even numbers $\mathbb{E}$ (the set of odd numbers  $\mathbb{O}$)  &  countable, $\aleph_0$ & $\frac{\G1}{2}$    \\
 &   &     \vspace{-2mm}   \\
 the set of integers $\mathbb{Z}$    &  countable, $\aleph_0$ & 2\G1+1   \\
 &   &     \vspace{-2mm}\\
    $\mathbb{Z} \setminus \{ 0 \} $  &  countable, $\aleph_0$ & 2\G1    \\
 &   &     \vspace{-2mm}\\
  squares of natural numbers $\mathbb{G} = \{ x : x= n^2, x \in \mathbb{N},\,\, n   \in \mathbb{N} \}$   &  countable, $\aleph_0$ & $ \lfloor \sqrt{\G1} \rfloor$      \\
 &   &     \vspace{-2mm}   \\
   pairs of natural numbers $\mathbb{P}  = \{ (p,q) : p   \in
\mathbb{N},\,\, q \in \mathbb{N} \}$   &  countable, $\aleph_0$ &  $\G1^2$   \\
   &   &     \vspace{-2mm}\\
 the set of numerals  $\mathbb{Q}_1  = \{ \frac{p}{q}:   p   \in \mathbb{Z}, \,\, q \in \mathbb{Z},\,\,\, q \neq
 0
 \} $   &  countable, $\aleph_0$ &  $4\G1^2+2\G1$   \\
  &   &     \vspace{-2mm}\\
 the set of numerals  $\mathbb{Q}_2  = \{ 0,   -\frac{p}{q}, \,\, \frac{p}{q} : p   \in
\mathbb{N}, \,\,q \in \mathbb{N} \} $   &  countable, $\aleph_0$ &  $2\G1^2+1$   \\
  &   &     \vspace{-2mm}\\
  the power set of the set
  of natural numbers  $\mathbb{N}$   &  continuum, \textfrak{c}  & $2^{\mbox{\scriptsize{\ding{172}}}}$     \\
 &   &     \vspace{-2mm}   \\
 the power set of the set
  of even numbers  $\mathbb{E}$   &  continuum, \textfrak{c} & $2^{0.5\mbox{\scriptsize{\ding{172}}}}$     \\
 &   &     \vspace{-2mm}   \\
 the power set of the set
  of integers  $\mathbb{Z}$   &  continuum, \textfrak{c} & $2^{2\mbox{\scriptsize{\ding{172}}}+1}$     \\
 &   &     \vspace{-2mm}   \\
the power set of the set
   of numerals  $\mathbb{Q}_1$    &  continuum, \textfrak{c} & $2^{\mbox{\scriptsize{$4\G1^2+2\G1$}}}$     \\
 &   &     \vspace{-2mm}   \\
the power set of the set
   of numerals  $\mathbb{Q}_2$    &  continuum, \textfrak{c} & $2^{\mbox{\scriptsize{$2\G1^2+ 1$}}}$     \\
 &   &     \vspace{-2mm}   \\
  numbers $x \in [0,1)$ expressible
in the binary   numeral system  &  continuum, \textfrak{c} & $2^{\mbox{\scriptsize{\ding{172}}}}$     \\
 &   &     \vspace{-2mm}   \\
 numbers $x \in [0,1]$ expressible
in the binary   numeral system  &  continuum, \textfrak{c} & $2^{\mbox{\scriptsize{\ding{172}}}}+1$     \\
 &   &     \vspace{-2mm}   \\
numbers $x \in (0,1)$ expressible
in the decimal   numeral system   &  continuum, \textfrak{c} & $10^{\mbox{\scriptsize{\ding{172}}}}-1$    \\
 &   &     \vspace{-2mm}   \\
numbers $x \in [0,2)$ expressible
in the decimal   numeral system  &  continuum, \textfrak{c} & $2 \cdot 10^{\mbox{\scriptsize{\ding{172}}}}$    \\
\hline
\end{tabular}
\end{center}
\end{table}

Let us consider now Table~\ref{table1} taken from \cite{EMS} that will be very useful in our study. It    shows that the \G1-based methodology allows one to count the number of elements of certain infinite sets with the precision of one element.

In order to proceed, let us informally define the set, $\mathbb{N}$, of natural numbers
\beq
\mathbb{N} = \{1, 2, 3, 4,  5, \,\, \ldots \,\, \}
\label{paradoxes_N}
 \eeq
as the set of numbers used to count objects.  Notice that nowadays not only positive integers are taken as elements of $\mathbb{N}$, but also zero is frequently included in $\mathbb{N}$ (see, e.g., \cite{Bagaria_1,Heller,Mancosu_2016}). However, since historically zero has been invented significantly later with respect to positive integers used for counting objects, zero is not included in
$\mathbb{N}$ in this article.

As it can be seen from  Table~\ref{table1}, the \G1-based methodology provides us with these results: the set  $\mathbb{N}$ has \G1 elements and the set of integers,   $\mathbb{Z}$, has 2\G1+1 elements. Thus, thanks to this novel way of counting,  it becomes possible to compute (and distinguish) the exact number of   elements of these two countable sets.  Notice (see \cite{EMS}, p. 241) that in this methodology \G1 is the last natural number and positive integers greater than grossone are called \emph{extended natural numbers}.  Another result from \cite{EMS} (once again see p. 241) that is important for our further consideration  regards the number of elements of the set,~$C_m$, of $m$-tuples of natural numbers:
 \beq
 C_m  =
\{
  (a_1, a_2, \ldots, a_{m-1}, a_m ) : a_i \in   \mathbb{N}, 1 \le i \le m
 \}, \hspace{3mm} 2 \le m \le
\G1.
 \label{tuples_2}
 \eeq
It is
known from combinatorial calculus that if we have $m$ positions and
each of them can be filled in by one of $l$ symbols, the number of
the obtained $m$-tuples is equal to $l^m$. In our case, since
$\mathbb{N}$  has grossone elements,~$l = \G1$. As a consequence, the set $C_m$
has $\G1^m$ elements. In particular, in the case   $m=\G1$,  the corresponding set $C_{\tiny{\G1}}$
has $\G1^{\tiny{\G1}}$ elements.

In order to conclude this brief tour in the \G1-based methodology let us show how arithmetical operations can be executed with infinite and infinitesimal numbers involving grossone. These numbers can have different infinite parts corresponding in their simplest form to finite positive powers of grossone. The \G1-based numbers can also have infinitesimal parts corresponding in their simplest form to finite negative powers of grossone. Finite  numbers $a$ are represented in the form $a=a\cdot \G1^0$ using the fact that $1=  \G1^0$ (see \cite{EMS} for a detailed discussion). Let us consider  as an example the following five numbers:
\[
\begin{array}{rcl}
A & = & 74.9\G1^{42.3}+5.1\G1^{0}+13.8\G1^{-25.6},\\
B & = &  5.7\G1^{16.8} -7.4\G1^{-14.9},\\
C & = &     74.9\G1^{42.3}+5.7\G1^{16.8}+5.1\G1^{0}-7.4\G1^{-14.9}+13.8\G1^{-25.6},\\
D & = &    74.9\G1^{42.3}-5.7\G1^{16.8}+5.1\G1^{0}+7.4\G1^{-14.9}+13.8\G1^{-25.6},\\
E & = & 426.93\G1^{59.1}-554.26\G1^{27.4}+29.07\G1^{16.8}+78.66\G1^{-8.8}\\
& & -37.74\G1^{-14.9}-102.12\G1^{-40.5}.
\end{array}
\]
The first of them, $A$, has one infinite part,  $74.9\G1^{42.3}$, one finite part, $5.1\G1^{0}$, and one infinitesimal part, $13.8\G1^{-25.6}$. The second number, $B$, has one infinite part, $5.7\G1^{16.8}$, and one infinitesimal part, $ -7.4\G1^{-14.9}$. The third number, $B$, has two infinite parts, one finite part, and two infinitesimal parts, etc. The arithmetic with \G1-based numbers works in such way (see \cite{EMS} for a formal detailed description of operations) that
\[
A + B = C,\hspace{3mm} A - B = D,\hspace{3mm} A \cdot B = E,\hspace{3mm} E / A = B.
\]

\section{\G1-based estimates}

We are ready now to use the \G1-based machinery to estimate the number of elements of the set, $A$, of algebraic numbers. We shall follow considerations made in the Introduction by substituting (where it is possible)  the cardinal number~$\aleph_0$ with the \G1-based quantities. The following theorem holds.

\begin{theorem}
The number of elements, $\hat{A}$, of  the set of all algebraic numbers  $A$ can be estimated as follows
 \beq
 4\G1+1 < \hat{A} <  \frac{(2\G1+1)((2\G1+1)^{\tiny{\G1}}(2\G1^2-1)+1)}{2\G1}.
 \label{estimate}
 \eeq
\end{theorem}

\textbf{Proof.}  We start by discussing the left-hand estimate in (\ref{estimate}). Since the set of integers, $\mathbb{Z}$, has 2\G1+1 elements, real roots of the equations $z-n=0,$ $n \in \mathbb{Z},$ give us the first 2\G1+1 algebraic numbers. Then,  equations $z^2-n=0$ for negative $n \in \mathbb{Z}$ give us other 2\G1 roots of the form $\pm \sqrt{n}= \pm \sqrt{|n|}i, i=\sqrt{-1},$ that do not coincide with the previously counted real roots. Thus, we have at least 4\G1+1 different algebraic numbers.

Let us discuss now the right-hand estimate in (\ref{estimate}) and recall once again that the set of integers, $\mathbb{Z}$, has 2\G1+1 elements. Then, since in (\ref{tuples_1}) numbers $a_i \in \mathbb{Z},  i=0,\ldots, k,$ with $a_0 \neq 0$, the first position in the tuple (\ref{tuples_1}) can be filled in by 2\G1 elements since $a_0 \neq 0$ and the remaining $k$ positions can be filled in by any from $2\G1+1$ integers. As a result, it follows from considerations related to  (\ref{tuples_2}) that the number of elements of the set $B_k$ of tuples (\ref{tuples_1}) has $2\G1(2\G1+1)^{k}$ elements.

The polynomial (\ref{polyn}) corresponding to each tuple (\ref{tuples_1}) has $k$ complex roots counted according to their multiplicities, i.e., the number, $l_k$, of different algebraic numbers   corresponding to this polynomial is    $l_k \le k$. Notice also that polynomials obtained from different tuples can have some common roots, as well.  Thus, the quantity of algebraic numbers,~$A_k$, corresponding to the set $B_k$ can be estimated as $A_k < 2\G1(2\G1+1)^{k}k$. In order to estimate the number, $\hat{A}$, of all algebraic numbers we should consider the union of the roots of polynomials linked to  all sets $B_k$. Recall now that~$k$ is a natural number and the set of natural numbers has~\G1 elements (see Table~\ref{table1}). Thus, it follows that $1 \le k \le \G1$ and we obtain the following estimate
\beq
\hat{A} < \sum_{k=1}^{\G1} A_k < \sum_{k=1}^{\G1}2\G1(2\G1+1)^{k}k = 2\G1\sum_{k=1}^{\G1}(2\G1+1)^{k}k.
 \label{estimate1}
 \eeq
Let us calculate the summation in (\ref{estimate1}) indicating it as
\[
S(\G1)=\sum_{k=1}^{\G1}(2\G1+1)^{k}k.
\]
Obviously, the sum $S(\G1)$ can be rewritten as
\[
S(\G1)=  (2\G1+1) + 2(2\G1+1)^2 + 3(2\G1+1)^3 + \ldots + (\G1-1)(2\G1+1)^{\tiny{\G1}-1} + \G1(2\G1+1)^{\tiny{\G1}}.
\]
Let us now divide this identity by 2\G1+1. We obtain
\[
\frac{S(\G1)}{2\G1+1}=     1+ 2(2\G1+1)  + 3(2\G1+1)^2 + \ldots + (\G1-1)(2\G1+1)^{\tiny{\G1}-2} + \G1(2\G1+1)^{\tiny{\G1-1}}
\]
and subtract the former relation from tha latter. We get
\[
\frac{S(\G1)}{2\G1+1}- S(\G1) =     1+  (2\G1+1)  +  (2\G1+1)^2 + (2\G1+1)^3 + \ldots
\]
\beq
  \ldots + (2\G1+1)^{\tiny{\G1}-2} +  (2\G1+1)^{\tiny{\G1-1}} - \G1(2\G1+1)^{\tiny{\G1}}.
\label{sum}
 \eeq
Positive summands in  (\ref{sum}) represent a geometric progression with the common ratio $q=2\G1+1$ (see \cite{EMS}, section 6.2, for a detailed discussion on summation with infinite and infinitesimal \G1-based numbers). As a result,  (\ref{sum}) can be rewritten as
\[
 \frac{S(\G1)}{2\G1+1} - S(\G1) =  \frac{1-(2\G1+1)^{\tiny{\G1}}}{1-(2\G1+1)}   - \G1(2\G1+1)^{\tiny{\G1}}
\]
from where we obtain
\[
-\frac{2\G1}{2\G1+1} S(\G1)  =  \frac{(2\G1+1)^{\tiny{\G1}}(1-2\G1^2)-1}{2\G1}.
\]
This relation allows us to express $ S(\G1)$ as follows
\[
  S(\G1)  =  \frac{(2\G1+1)((2\G1+1)^{\tiny{\G1}}(2\G1^2-1)+1)}{4\G1^2}.
\]
It is sufficient now to substitute this result in (\ref{estimate1}) to obtain the right-hand estimate in (\ref{estimate}). This fact concludes the proof. \hfill $\Box$

\textbf{Remark.} Let us emphasize that the lower estimate in (\ref{estimate})   can be further improved by   counting  algebraic numbers being roots of other polynomials.  For example, let us consider polynomials $z^2-n=0$, where $n$ is one of the first $k$ positive prime numbers
\[
\underbrace{2, 3, 5, 7, 11, 13, 17, 19, 23, 29, 31, 37, 41, \ldots }_{ k\,\,  \mbox{{\scriptsize first primes}}}
\]
Then, clearly, the $2k$ roots of these $k$ polynomials are not among $4\G1+1$ numbers counted in the theorem (e.g.,  numbers $\pm\sqrt{2}$   being roots of the polynomial  $z^2-2=0$   have not been counted). Thus,  the lower estimate in (\ref{estimate}) can be re-written as $4\G1+2k+1$.

\section{A brief conclusion}
With respect to  the set   of algebraic numbers, $A$, the classical terminology of cardinals proposed by Cantor says that $A$  is countable. This result is obtained by using the fact that a set being a countable union of  countable sets  is also countable. In the present paper, instead of the usage of Cantor's cardinals,   a recently introduced   methodology working with \G1-based infinite numerals is applied. This has been done because the \G1-based  methodology   can provide  a more accurate measurement  of certain infinite sets (see Table~\ref{table1}).   For example,  if two infinite countable sets $P$ and $P'$ are considered such that   the set $P'$ has been obtained by adding to $P$ one element $p \notin P$ then cardinals of Cantor allow us to say only that both sets are countable. In fact, it follows that the cardinality of $P$ is $\aleph_0$ and the cardinality of $P'$  is $\aleph_0+1=\aleph_0$. In other words, the numeral system of Cantor does not allow us to register the fact that one element has been added to $P$ to obtain $P'$. In contrast, if we consider the set $P$ and its number of elements $n$ is expressed in  \G1-based numerals, then the set $P'$   has $n+1>n$ elements and \G1-based  numerals allow us to register this fact (see, for example, the fourth and the fifth lines in Table~\ref{table1}). Thus, we have a complete  analogy  with what   happens with finite sets: when one adds to a finite set $F$ having $n$ elements an element    $f \notin F$ then the resulting set $F'$   has $n+1>n$ elements.

Due to this more precise analysis of infinite sets,   \G1-based  numerals together with some combinatorial considerations have allowed us to obtain lower and, what is more important, upper bounds for the number of elements of the set of algebraic numbers $A$. In general, it can be safely said that  considerations involving grossone-based numerals performed in this paper open   new promising prospectives  for measuring both countable sets and sets that are countable  unions   of countable sets.

\section*{Declarations}
The author states that there is no conflict of interest and no funding has been used to execute this research.

\section*{Data availability statement} Data sharing is not applicable to this article as no datasets were generated or analysed during the current study.





\end{document}